\documentclass{elegantpaper}

\usepackage{amsmath}
\usepackage{amsfonts}
\usepackage{amssymb}
\usepackage{mathrsfs}
\usepackage{graphicx}
\usepackage{xcolor}

\newcommand  \bglb {\big (}
\newcommand  \bgrb {\big )}
\newcommand  \bglc {\big \{}
\newcommand  \bgrc {\big \}}

\newcommand   \Integers {\mathbb Z}
\newcommand   \PositiveIntegers {{\mathbb Z}^{+}}

\newcommand  \xx   {\mathbf x}
\newcommand  \yy   {\mathbf y}

\newcommand  \basel[2]{#1_{_{#2}}}

\newcommand  \tab  {\hspace*{0.5cm}}

\definecolor{titlecolor}{RGB}{144,48,48}
\definecolor{authorname}{RGB}{16,96,16}%
\definecolor{addrscolor}{RGB}{60,113,183}
\definecolor{secheader}{RGB}{16,32,128}%
\definecolor{refscolor}{RGB}{16,64,128}%

\begin{document}

\title{\textcolor{titlecolor}{\bf{Compactness of any Countable Product of Compact Metric Spaces in Product Topology without Using Tychonoff’s Theorem}}}

\author{\\
\\
{\textcolor{authorname}{\Large{\bf{Garimella Sagar$^{1}$}}}} \tab
 {\textcolor{authorname}{and}} \tab
{\textcolor{authorname}{\Large{\bf{Duggirala Ravi$^{2}$}}}}\\
\\
 {\textcolor{addrscolor}{\scriptsize{\bf{$^{1}$E-mail ~: \tab Sagar.Garimella@gmail.com}}}} \\
 {\textcolor{addrscolor}{\scriptsize{\bf{$^{2}$E-mail ~: \tab ravi@gvpce.ac.in; ~~ drdravi2000@yahoo.com;
 ~~	 duggirala.ravi@rediffmail.com; ~~duggirala.ravi@yahoo.com} }}}
\\
\\
}

\date{}

\maketitle

\textcolor{secheader}{
\begin{abstract}
For infinite products of compact spaces, Tychonoff’s theorem asserts that their product is compact, in the product topology. Tychonoff’s theorem is shown to be equivalent to the axiom of choice. In this paper, we show that any countable product of compact metric spaces is compact, without using Tychonoff’s theorem. The proof needs only basic and standard facts of compact metric spaces and the Bolzano-Weierstrass property. Moreover, the component spaces need not be assumed to be copies of the same compact metric space, and each component space can be an arbitrary nonempty compact metric space independently.
Total boundedness together with completeness of a metric space implies its compactness. Completeness of a product of complete spaces is easily inferred from the completeness of each component. Total boundedness therefore suffices to prove the compactness of a product space consisting of countably (infinitely) many nonempty compact component spaces. The countable infiniteness is needed in the proof to exhibit a standard metric that gives rise to the product topology. Any such metric topology for the product arises as exhibited, and they are all equivalent to the product topology. The requirement of summability of the sequences restricts the scope of the result to countably infinite products.  In summary, the product space obtained by taking the product of any sequence of nonempty compact metric spaces in the product topology is shown to be compact, using only the basic and standard facts of compact metric spaces. In conclusion, compactness of the product of a countably infinitely many nonempty compact metric spaces can be proved within Cantor’s set theory, without using the axiom of choice and Tychonoff’s theorem.
\end{abstract}
}

\begin{small}
\textcolor{secheader}{
}
\end{small}

\textcolor{titlecolor}{\section{\bf{Introduction}}}
Tychonoff’s theorem states that any infinite product of compact spaces is compact in the product topology. The only nonessential assumption in the general statement is that the component spaces are copies of the same compact space. It is interesting to explore the possibility of proving compactness without using Tychonoff’s theorem, for infinite products.  This is possible, but with some restrictions. It is assumed that the product consists of only countably many nonempty compact metric spaces.  For metric spaces, compactness is equivalent to the Bolzano-Weierstrass property, whereby requiring that any sequence containing infinitely many distinct elements has a limit point. By completeness, the limit point of any such sequence is required to belong to the given space. Metrizability of a countable product of nonempty compact metric spaces is easy to establish. The completeness of the product of complete metric spaces can be inferred by the componentwise completeness, for Cauchy sequences in metric spaces. It follows that the total boundedness implies compactness, for a complete metric space. Now,  any $\epsilon$-net cover, for the product, induces an $\frac{\epsilon}{2}$-net cover for a product of a fixed finite number of component spaces, which permits a finite $\frac{\epsilon}{2}$-net cover, as the components are compact themselves, enabling us to recover a finite $\epsilon$-net cover of the originally given product space consisting of countably infinitely many nonempty compact metric spaces, and hence the product space itself becomes totally bounded. In this regards, the main fact is the demonstration of a metric that is easy to define on the product space, but restricting its applicability to sequences of spaces, meaning countably infinitely many nonempty compact metric spaces. The metric topology thus obtained is equivalent to the product topology.

\textcolor{titlecolor}{\section{\label{Sec-Main-Result}\bf{Main Result}}}
Let $\Integers$ and $\PositiveIntegers$ be the sets of integers and positive integers, respectively.
For the notation, we follow mostly the description given in \cite{Sutherland-2009}.
A topological space is denoted by $(X,~\basel{\mathcal T}{X})$, where $X$ is a nonempty space
and $\basel{\mathcal T}{X}$ is a topology on $X$. For a metric space $X$, with metric $\basel{d}{X}$,
the notation $(X, \, \basel{d}{X})$ is reserved. The following facts are needed later:
\begin{enumerate}
\item For a metric space, the compactness is equivalent to its satisfying the Bolzano-Weierstass propoerty \cite{Croom-2008}.
\item For $\bglb\basel{X}{\alpha}, ~ \basel{d}{\basel{X}{\alpha}}\bgrb$, $\alpha \in \Lambda$, for some index set $\Lambda$,
let ${\mathbf X} = \prod_{\alpha \in \Lambda} \basel{X}{\alpha}$, equipped with the product topology. If $X$ is metrizable,
such that the projection maps onto the component spaces remain continuous, and $\basel{X}{\alpha}$ is a complete metric space,  for every $\alpha \in \Lambda$, then $X$ is also a complete metric space.
\end{enumerate}
The proof of the second item is trivial, as the projection maps onto the component spaces must remain continuous,
and Cauchy sequences in the product space ${\mathbf X}$ are mapped into Cauchy sequences for each component space.
The component spaces being complete, the uniqueness of a limit point of the Cauchy sequence ensures that  the components
of the limit of the sequence in ${\mathbf X}$ must be the same as given by the limit points in the component spaces.
The subscript $X$ in $\basel{\mathcal T}{X}$ and $\basel{d}{X}$ are sometimes omitted for the sake of simplicity of notation.

\begin{proposition}
For any metric space $(X, \, d)$, there is a bounded metric $d'$ on $X$,
generating the same topology.
\end{proposition}
\proof The proof is based on an exercise given in \cite{Lang-RA-1983}.
Let $d'$ be defined by
\[
d'(x,\, y) = \left \{
\begin{array}{lcl}
d(x, \, y) &,&  \textrm{if}~~  d(x, \, y) < 1\,, ~~\textrm{and}\\
1 &,&  \textrm{if}~~  d(x, \, y) \geq  1
\end{array} \right.
\]
for $x, y \in X$. Then, $d'(x,\, y) = 0$, if and only if $d(x,\, y) = 0$,  and $d'(x,\, y) = d'(y,\, x)$,
for $x, y \in X$. For the triangle inequality, let $x,\, y,\,  z \in X$, and the following
two possibilities are considered:
\begin{enumerate}
\item if $d'(x,\, y) = 1$ or $d'(y,\, z) = 1$, then $d'(x,\, z) \leq 1 \leq d'(x,\, y) + d'(y,\, z)$, and

\item if $d'(x,\, y) < 1$ and $d'(y,\, z) < 1$, then $d'(x,\, y) = d(x,\, y)$ and $d'(y,\, z) = d(y, \, z)$, and
$d'(x,\, z) \leq d(x,\, z) \leq d(x,\, y) + d(y,\, z) = d'(x,\, y) + d'(y,\, z)$.
\end{enumerate}
The basic open sets of radius less than $1$ are the same with respect to the metrics $d$ and $d'$. \qed
\\
 In the next couple of propositions, a more general method for
 construction of a bounded metric, for  a given metric is described.
\begin{proposition}
Let $(X,\, d)$ be a metric space, and $h\,:\,[0,\, \infty) \to [0,\, a)$
be a  nondecreasing continuous function, for some fixed real number $a > 0$,
 such that $h(0) = 0$, $h(b+c) \leq h(b)+h(c)$, for every nonnegative
 real numbers $b$ and $c$, and $h$ is one-to-one and strictly increasing
 when restricted to the interval $[0,\, \eta]$, for some fixed real number
 $\eta > 0$.  Let $d'(x,\, y) = h\bglb (d(x ,\, y)\bgrb$,  for $x,\, y \in X$.
  Then, $d'$ is a metric on $X$, and the metric spaces $(X, \, d)$ and $(X, \, d')$
are topologically equivalent.
\end{proposition}
\proof Since $h$ is nondecreasing and strictly increasing on $[0,\, \eta]$, it follows
 that, for any nonnegative real number $b$,  $h(b) = 0$ if and only if $b = 0$.
Symmetry of $d'$ follows from that of $d$, and for the  triangle inequality,
it is  observed that, for $x,\, y, \, z \in X$,
\begin{eqnarray*}
d'(x,\,  z)  & = & h\bglb d(x,\, z) \bgrb \\
           & \leq  & h\bglb d(x,\, y) +  d(y,\, z)  \bgrb\,, \tab \textrm{since} ~ h ~ \textrm{is nondecreasing} \\
          & \leq  & h\bglb d(x,\, y) \bgrb + h\bglb d(y,\, z)  \bgrb\,, \tab \textrm{by the special condition on} ~ h ~ 
\end{eqnarray*}
 The basic open sets of radii $\epsilon$ for the metric $d$ and
  $h (\epsilon)$ for the metric $d'$, for every real  number
 $ \epsilon$, such that  $0 < \epsilon < \eta $, are identical. \qed
 \\

The next proposition provides a means for constructing
a diffeomorphism of the metrics $d$ and $d'$:
\begin{proposition}
With the  notation of the preceding proposition, if the  homeomorphism $h$ of the metrics
$d$ and $d'$ satisfies the condition that its derivative $h'$ is defined almost everywhere
and nonincreasing on $[0 ,\, \infty)$, then the condition  $h(b+c) \leq h(b)+h(c)$,
 for every nonnegative real numbers $b$ and $c$, holds.
\end{proposition}
 \proof Since  $h(0) = 0$  and $h'$ is nonincreasing, $h'(x+c) \leq h'(x)$,
 {\em a.e.} for $0\leq x \leq b$, it holds that
  $\int_{0}^{b} h'(x+c) dx \leq   \int_{0}^{b} h'(x) dx$, and that
   $h(b+c)-h(c) \leq h(b)$, whereby the required condition follows.  \qed
\paragraph{\underline{Attention}.}~ Since $h$ is nondecreasing, $h' \geq 0$, {\em a.e.} on $[0,\, \infty)$.
\\

A metric space $(X,\, d)$ is said to be totally bounded, if for every $\epsilon > 0$,
there is finite subset $A(\epsilon)\subseteq X$, such that
 $X \subseteq \bigcup_{x \in A(\epsilon)} B(x, \, \epsilon)$, where 
 $B(x, \, \epsilon) = \bglc y \in X \, :\, d(x, y) < \epsilon\bgrc$. 
A complete metric space is compact if and only if it is totally bounded:
\begin{proposition}
Let $(X,\, d)$ be a complete metric space. If $X$ is totally bounded, then $X$ is compact.
Conversely, if $(X,\, d)$ is a compact metric space, then $X$ is complete and totally bounded.
\end{proposition}
\proof Following \cite{Croom-2008}, the compactness is proved by showing that the 
Bolzano-Weierstrass property holds for $X$. By the assumption of total boundedness,
for every positive integer $n \in \PositiveIntegers$, there is a finite subset
$\basel{A}{n} \subseteq X$, such that
$X \subseteq \bigcap_{y \in \basel{A}{n}} B\bglb y,\, \frac{1}{2n}\bgrb$.
Let $\{ \basel{x}{i} \,:\, i \in \PositiveIntegers\}$ be a sequence, consisting
of infinitely many distinct elements.  By the finiteness of $\basel{A}{n}$,
there is a $\basel{y}{n} \in \basel{A}{n}$, such that
 $B\bglb\basel{y}{n},\, \frac{1}{2n}\bgrb$ contains infinitely many elements
 of the given sequence. For any fixed $n \in \PositiveIntegers$,
 another element $\basel{y}{n+1} \in \basel{A}{n+1}$ is obtained, such that 
  $B\bglb\basel{y}{n+1},\, \frac{1}{2(n+1)}\bgrb \cap B\bglb\basel{y}{n},\, \frac{1}{2n}\bgrb$
  contains infinitely many elements from the subsequence that is already included in
$B\bglb\basel{y}{n},\, \frac{1}{2n}\bgrb$. Now,
$d \bglb\basel{y}{n+i}, \, \basel{y}{n}\bgrb < \frac{1}{2n} + \frac{1}{2(n+i)} < \frac{1}{n}$,
for every $i \in \PositiveIntegers$, and the sequence $\basel{y}{n}$, $n \in \PositiveIntegers$,
is a Cauchy sequence. Let $z$ be the limit point of $\basel{y}{n}$, $n \in \PositiveIntegers$,
which is guaranteed to exist by the completeness of $X$. It can be easily checked that $z$
is also a limit point of the given sequence $\basel{x}{i}$, $i \in \PositiveIntegers$. 

   Conversely, if $(X, \,   d)$ is a compact metric space,
   then it satisfies the Bolzano-Weierstrass property  \cite{Croom-2008}, and hence,
   must be complete as a metric space. For the total boundedness, it is observed that
   $\bglc B(x, \, \epsilon)  \, : \,   x \in X \bgrc$ is an open cover of $X$, that
   admits a finite subcover, for every $\epsilon > 0$.  \qed
\\

Let$(\basel{X}{i}, \, \basel{d}{i})$, $i \in \PositiveIntegers$, be nonempty metric spaces,
with bounded metrics, such that for some positive integers $\basel{M}{i}$,
 $\basel{d}{i}(x, \, y) < \basel{M}{i}$, for all $x, \, y \in \basel{X}{i}$ and
$i \in \PositiveIntegers$.  Let $\basel{l}{i} > 0$ be positive real numbers such that
$L = \sum_{i = 1}^{\infty} \basel{l}{i} <\infty$. Let
 ${\mathbf X} = \prod_{i = 1}^{\infty} \basel{X}{i}$, with its topology
 $\basel{\mathcal T}{\mathbf X}$  given by the product. 
 For $\xx = (\basel{x}{1}, \, \basel{x}{2}, \, \basel{x}{3},\, \ldots ),\,  \yy = (\basel{y}{1}, \, \basel{y}{2}, \, \basel{y}{3},\, \ldots ) \in {\mathbf X}$,  let
 $D(\xx, \, \yy) = \sum_{i = 1}^{\infty} \frac{\basel{l}{i} \basel{d}{i}(\basel{x}{i},\, \basel{y}{i})}{\basel{M}{i}}$.
Since $\basel{l}{i} > 0$, for every $i \in \PositiveIntegers$, it follows that 
$D(\xx, \, \yy) = 0$, if and only if $\basel{d}{i}(\basel{x}{i},\, \basel{y}{i}) = 0$,
 for every $i \in \PositiveIntegers$. Symmetry and triangle inequality are easily seen to hold,
 to realize that $D$ is a metric on ${\mathbf X}$. A more important fact is that the metric topology,
 with respect to $D$, is equivalent to the product topology $\basel{\mathcal T}{\mathbf X}$, as shown
 in the following:
 
 \begin{theorem}
Let $(\basel{X}{i}, \, \basel{d}{i})$, $i \in \PositiveIntegers$, be nonempty metric spaces,
with bounded metrics, {\em i.e.},  $\basel{d}{i}(x, \, y) < \basel{M}{i}$, for $x, \, y \in \basel{X}{i}$ and
$({\mathbf X}, \, D)$ be the the metric space obtained as just described.
 Then, the topological space $({\mathbf X}, \, \basel{\mathcal T}{\mathbf X})$ and
the metric space $({\mathbf X}, \, D)$ are topologically equivalent.
\end{theorem}
\proof It may be recalled that $\basel{\mathcal T}{\mathbf X}$ is the product topology
on ${\mathbf X}$, with respect to projections onto the component spaces, while  $({\mathbf X}, \, D)$
is  the metric space with respect to the metric
 $D(\xx, \, \yy) = \sum_{i = 1}^{\infty} \frac{\basel{l}{i} \basel{d}{i}(\basel{x}{i},\, \basel{y}{i})}{\basel{M}{i}}$.
 for $\xx = (\basel{x}{1}, \, \basel{x}{2}, \, \basel{x}{3},\, \ldots ),\,  \yy = (\basel{y}{1}, \, \basel{y}{2}, \, \basel{y}{3},\, \ldots ) \in {\mathbf X}$, and for some sequence of positive real numbers  $\basel{l}{i} > 0$, $i \in \PositiveIntegers$.

 Let $V \in \basel{\mathcal T}{\mathbf X}$ be a nonempty open set with respect to the product topology,
and $\xx \in V$ be an arbitrary element. By the condition required by the product topology, there is a positive
integer $n \in \PositiveIntegers$, such that $\basel{U}{i} = \basel{X}{i}$, whenever $i \geq n+1$.
Let $\epsilon > 0$ be such that, for any sequence $\yy \in {\mathbf X}$, if
 $\frac{\basel{l}{i} \basel{d}{i}(\basel{x}{i},\, \basel{y}{i})}{\basel{M}{i}} < \epsilon$,
for $1 \leq i \leq n$,  then $\yy \in V$.  Let $B(\xx, \, \epsilon)$ be the basic open subset
 $\{\yy \in {\mathbf X} \, :\, D(\xx,\, \yy) < \epsilon\}$. Then, $B(\xx, \, \epsilon) \subseteq V$,
 and the latter is an open set with respect to the metric topology.
 For the converse inclusion, let $B(\xx, \, \epsilon)$ be a basic open subset of ${\mathbf X}$,
 with respect to the metric topology, for some $ \xx \in {\mathbf X}$ and $\epsilon > 0$.
 Let $n \in \PositiveIntegers$ be sufficiently large,
 such that $\sum_{i = n+1}^{\infty} \basel{l}{i} < \frac{\epsilon}{2}$.
 Then, the open  set $V = \prod_{i = 1}^{\infty} \basel{U}{i}$, defined by
 constraining the first $n$ components $\basel{y}{i}$, $1 \leq i \leq n$, by the condition 
$\sum_{i = 1}^{n} \frac{\basel{l}{i} \basel{d}{i}(\basel{x}{i},\, \basel{y}{i})}{\basel{M}{i}} < \frac{\epsilon}{2}$,
and letting $\basel{U}{i} = \basel{X}{i}$, for $i \geq n+1$, is included in $B(\xx, \, \epsilon)$,
and the latter is open with respect to the product topology.  \qed
\paragraph{\underline{Attention}.}~~ The metric spaces $\bglb \basel{X}{i},\, \basel{d}{i}\bgrb$
are not necessarily assumed to be complete or  compact, for any index $i\in \PositiveIntegers$, but
on that $\basel{d}{i}(x,\, y) \leq \basel{M}{i}$, for some real numbers $\basel{M}{i} > 0$,
 for $i \in \PositiveIntegers$. The metric topology of $\bglb {\mathbf X}, \, D\bgrb$ is still
 equivalent to the product topology, where ${\mathbf X} = \prod_{i = 1}^{\infty} \basel{X}{i}$ and
 $D(\xx, \, \yy) = \sum_{i = 1}^{\infty} \frac{\basel{l}{i} \basel{d}{i}(\basel{x}{i},\, \basel{y}{i})}{\basel{M}{i}}$.
 for $\xx = (\basel{x}{1}, \, \basel{x}{2}, \, \basel{x}{3},\, \ldots ),\,  \yy = (\basel{y}{1}, \, \basel{y}{2}, \, \basel{y}{3},\, \ldots ) \in {\mathbf X}$, and for some sequence of positive real numbers  $\basel{l}{i} > 0$, $i \in \PositiveIntegers$,
\\

 \begin{theorem}
 {\bf {(Sagar)}}~~
 Let$(\basel{X}{i}, \, \basel{d}{i})$, $i \in \PositiveIntegers$, be nonempty compact metric spaces,
 and $({\mathbf X}, \, D)$ be the the metric space obtained as just described.
 The metric space $({\mathbf X}, \, D)$ is totally bounded,
 and becomes compact, with respect to the product topology.      
\end{theorem}
\proof For any $\epsilon > 0$, the objective is to show that
there is a finite set ${\mathcal A}(\epsilon) \subseteq {\mathbf X}$, such that
${\mathbf X} \subseteq \bigcup_{\xx \in {\mathcal A}(\epsilon)} B(\xx, \, \epsilon)$.
As in the converse inclusion  part of the proof of
the preceding theorem, let $n \in \PositiveIntegers$ be sufficiently large,
 such that $\sum_{i = n+1}^{\infty} \basel{l}{i} < \frac{\epsilon}{2}$.
 Then, the open  set $V = \prod_{i = 1}^{\infty} \basel{U}{i}$, defined by
 constraining the first $n$ components $\basel{y}{i}$, $1 \leq i \leq n$, by the condition 
$\sum_{i = 1}^{n} \frac{\basel{l}{i} \basel{d}{i}(\basel{x}{i},\, \basel{y}{i})}{\basel{M}{i}} < \frac{\epsilon}{2}$,
and letting $\basel{U}{i} = \basel{X}{i}$, for $i \geq n+1$, is included in $B(\xx, \, \epsilon)$.
By the compactness of $\prod_{i = 1}^{n} \basel{X}{i}$, that can be proved using only standard facts,
without using the Tychonoff's theorem, there is a finite set $\basel{A}{n}\bglb \frac{\epsilon}{2} \bgrb$,
such that for every element  $(\basel{y}{1}, \, \ldots, \, \basel{y}{n}) \in \prod_{i = 1}^{n} \basel{X}{i}$,
there an element $(\basel{x}{1}, \, \ldots, \, \basel{x}{n}) \in \basel{A}{n}\bglb \frac{\epsilon}{2}\bgrb$,
such that
$\sum_{i = 1}^{n} \frac{\basel{l}{i} \basel{d}{i}(\basel{x}{i},\, \basel{y}{i})}{\basel{M}{i}} < \frac{\epsilon}{2}$, 
and hence, every element in $ \prod_{i = 1}^{n} \basel{X}{i}$ is at a distance less  than $\frac{\epsilon}{2}$,
from or to an element in  $ \basel{A}{n}\bglb \frac{\epsilon}{2}\bgrb $, with respect to the metric 
$\sum_{i = 1}^{n} \frac{\basel{l}{i} \basel{d}{i}(\basel{x}{i},\, \basel{y}{i})}{\basel{M}{i}}$.
The collection ${\mathcal A}(\epsilon)$ of sequences in ${\mathbf X}$, obtained by taking
the first $n$ components to be those of the finite sequences in $\basel{A}{n}\bglb \frac{\epsilon}{2} \bgrb$
and the remaining components  to be some specified default elements from the corresponding component spaces,
can be shown to be a finite  $\epsilon$-net for ${\mathbf X}$. \qed

\textcolor{titlecolor}{
\section{\label{Sec-Apps-and-Misc}Applications and Miscellany} 
}
In this section, some applications of the results of the preceding section,
with reference to quotient topology, for development of understanding of
new spaces are presented.
\\

\textcolor{titlecolor}{
\subsection{\label{Sec-Bas-Top-Fac}Some Basic Facts of Topology} 
}

The following couple of propositions are standard  facts and will be needed later:

\begin{proposition}
Let $X$ and $U$ be topological spaces, and $f\,:\, X \to Y$ be a continuous function
from $X$ onto $Y$, with $f(X) = Y$. If $X$ is compact,  then so is $Y$.
If $f$ is one-to-one and  $Y$  is Hausdorff, then $X$ is also Hausdorff.
\end{proposition}
\proof For any open cover $\bglc \basel{U}{\alpha}\,:\, \alpha \in  \Lambda \bgrc$
of $Y$, for any index set $\Lambda$, $\bglc f^{-1}\bglb \basel{U}{\alpha}\bgrb\,:\, \alpha \in  \Lambda \bgrc$
is an open cover of $X$, by the continuity of $f$. There is a positive integer $n \in \PositiveIntegers$
and indexes $\basel{\alpha}{i}\in \Lambda$, for $1 \leq i \leq n$, such that
$X \subseteq f^{-1}\bglb\basel{U}{\basel{\alpha}{1}}\bgrb \cup \,  \ldots   \,\cup  f^{-1}\bglb\basel{U}{\basel{\alpha}{1}}\bgrb$,
by the compactness of  $X$.  Clearly, 
$Y \subseteq \basel{U}{\basel{\alpha}{1}}\cup \, \ldots \, \cup \basel{U}{\basel{\alpha}{n}}$,
and  $Y$ is compact.

  For the remaining part, let $\basel{x}{1},\,\basel{x}{2}\in X$, with $\basel{x}{1} \neq \basel{x}{2}$.
Since $f$ is one-to-one, it follows that $f\bglb\basel{x}{1}\bgrb \neq f\bglb\basel{x}{2}\bgrb$,
and since  $Y$ is  Hausdorff, there are open sets $\basel{U}{1},\,\basel{U}{2} \subseteq Y$,
with $f\bglb\basel{x}{1}\bgrb \in  \basel{U}{1}$ and $f\bglb\basel{x}{2}\bgrb \in  \basel{U}{2}$,
such that $\basel{U}{1} \cap \basel{U}{2}  = \emptyset$. Now, $f^{-1}\bglb\basel{U}{1}\bgrb$ and
$f^{-1}\bglb\basel{U}{1}\bgrb$ are open subsets of  $X$, with $\basel{x}{1} \in f^{-1}\bglb\basel{U}{1}\bgrb$
and  $\basel{x}{2} \in f^{-1}\bglb\basel{U}{2}\bgrb$, such that
 $f^{-1}\bglb \basel{U}{1}\bgrb \cap f^{-1}\bglb \basel{U}{2}\bgrb  = \emptyset$,
 and therefore, $X$ is Hausdorff.      \qed
\\

\begin{proposition}
Let $X$ and $U$ be topological spaces, and $f\,:\, X \to Y$ be a ine-to-one
continuous function from $X$ onto $Y$, with $f(X) = Y$. If $X$ is compact and
$U$ is  Hausdorff, then $f$ is a homeomorphism.
\end{proposition}
\proof It is required to prove that $f$ is an open map. Let $U \subseteq  X$
be an set. Then, $A = X \backslash U $ is closed, and, as
$X$ is compact,  $A$ is also compact. Now, $f(A)$ is compact,
and, since $Y$ is Hausdorff, it follows that $f(A)$ must be
closed \cite{Croom-2008}. As $f$ is a surjection,
the set $Y \backslash f(A) = f(X)\backslash f(A)$ 
is open. Since $f$ is one-to-one, it follows that
$f(X \backslash A) = f(X)\backslash f(A)$, and that
$f(U) = f(X \backslash A)$ is an open subset of $Y$.        \qed

\textcolor{titlecolor}{
\subsection{\label{Sec-Quot-Top}Quotient Topology} 
}

\begin{theorem}
Let $X$ and $Y$ be nonempty topological spaces, and $f\,:\, X \to Y$ be 
a continuous function, with $f(X0 = Y$. Let  $\equiv (~\textrm{w. r. t.} ~ f~)$
be the equivalence  relation on $X$, defined by the condition that,
for $\basel{x}{1},\, \basel{x}{2} \in X$, 
$\basel{x}{1} \equiv  \basel{x}{2} ~ (~\textrm{w. r. t.} ~ f~) $  if and only if
$f\bglb \basel{x}{1}\bgrb = f\bglb \basel{x}{2}\bgrb $. Let $W$ be the quotient space
$X / \equiv (~\textrm{w. r. t.} ~ f~)$, obtained by collecting the equivalence classes
with respect to $\equiv (~\textrm{w. r. t.} ~ f~)$, equipped with the quotient
topology, and let $\basel{\pi}{\equiv (~\textrm{w. r. t.} ~ f~)}$ be the
quotient map from $X$ onto $W$.  Let $g\,:\, W \to Y$ be the mapping defined by
$g(w) = f(x)$,  for any $w \in W$ and $x \in w$. The following statements hold:
\begin{enumerate}
\item the mapping $g$ is a well-defined one-to-one continuous function from $X$ onto $Y$\,;
\item if $X$ compact, then so are $Y$ and $W$\,;
\item if $Y$ is Hausdorff, then so is $W$\,;  ~~ and
\item if both $X$ is compact and $Y$ is Hausdorff, then $g$ is 
      a homeomprohism from  $X$ onto $W$.
\end{enumerate}
\proof The function $f$ is exactly $g\circ\basel{\pi}{\equiv (~\textrm{w. r. t.} ~ f~)}$,
hence $g$ is  continuous \cite{Sutherland-2009}. The remaining follow immediately
from the propositions of the preceding section.  \qed

\end{theorem}

\paragraph{\underline{Example}.}~~Let  $(\basel{X}{i},\, \basel{d}{i})$,  $i \in \PositiveIntegers$,
be the same copies of $(X,\, d)$, where $X = \{0,\,  1\}$ and $d(x,\,y) = |x-y|$, for $x,\, y \in X$. 
Let ${\mathbf X} = \prod_{i = 1}^{\infty}\basel{X}{i}$ and, for
$\xx = (\basel{x}{1},\,\basel {x}{2}, \,\basel{x}{3},\,\ldots ),\, 
\yy = (\basel{y}{1},\,\basel {y}{2}, \,\basel{y}{3},\,\ldots )  \in  {\mathbf X}$,
let $D(\xx,\, \yy) = \sum_{i = 1}^{\infty} 2^{-i} \basel{d}{i}(\basel{x}{i},\, \basel{y}{i})$.
Let $f\,:\, {\mathbf X} \to [0,\, 1]$,  be the mapping defined by
$f(\xx) =  \sum_{i = 1}^{\infty} 2^{-i}\basel{x}{i}$ 
for $\xx = (\basel{x}{1},\,\basel {x}{2}, \,\basel{x}{3},\,\ldots )  \in  {\mathbf X}$.
Now, ${\mathbf X}$ is compact, $[0,\,1]$ is Hausdorff, with respect to the standard absolute difference metric,
and  since  $|f(\xx) - f(\yy)| = $    
 $ \left | \sum_{i = 1}^{\infty} 2^{-i}\basel{x}{i} \, -  \,   \sum_{i = 1}^{\infty} 2^{-i}\basel{y}{i} \right |$
 $\leq $
$ \sum_{i = 1}^{\infty} 2^{-i}  |\basel{x}{i}  -   \basel{y}{i}|$,
it follows that $f$ is continuous.
Let $\equiv (~\textrm{w. r. t.} ~ f~)$ be the equivalence relation  on ${\mathbf X}$
as in the preceding theorem, and ${\mathbf W} = {\mathbf X} / \equiv (~\textrm{w. r. t.} ~ f~)$  be
the quotient  space, with the corresponding projection map  $\basel{\pi}{\equiv (~\textrm{w. r. t.} ~ f~)}$.
Then, the function $g\,:\, {\mathbf W} \to [0,\,1]$, defined as in the preceding theorem, is a homeomorphism.

\textcolor{titlecolor}{
\section{Conclusions} 
}
The compactness of the product of a sequence of compact metric spaces can be proved without using 
Tychonoff's theorem.   The diagonal sampling method presented in \cite{Rudin-1976} shows that 
Bolzano-Weierstrass property  holds for the product of a sequence of compact metric spaces.
In this context, the main contribution is the demonstration of a metric on the product of
a sequence of compact metric spaces, such that   the metric topology is equivalent to the
product topology.

 \textcolor{refscolor}{
}
\end{document}